\documentclass[12pt]{article}

\usepackage[a4paper, total={18cm, 24cm}]{geometry}

\usepackage[bottom]{footmisc}

\usepackage{mathtools}
\usepackage{amsmath,amssymb,amsbsy}
\usepackage{graphicx}
\usepackage{multirow}

\newtheorem{thm}{Theorem}[section]
\newtheorem{lem}[thm]{Lemma}
\newtheorem{obs}[thm]{Observation}

\newtheorem{cor}[thm]{Corollary}

\newtheorem{dfn}[thm]{Definition}

\newtheorem{conj}[thm]{Conjecture}

\makeatletter
\def\imod#1{\allowbreak\mkern10mu({\operator@font mod}\,\,#1)}
\makeatother

\newcommand{\zet}{\mathbb{Z}}
\def\gr{\mathcal{G}}
\newcommand{\gA}{\mathcal{A}}

\newcommand{\qed}{\hfill \rule{.1in}{.1in}}

\def\gcd{\mathop{\rm gcd}\nolimits}

\begin{document}
\title{Partition of Abelian groups into zero-sum sets by complete mappings and its application to the existence of a magic rectangle set}

\author{Sylwia Cichacz  \\
AGH University of Krak\'ow, Poland}

\maketitle
\begin{abstract}
A complete mapping  of a group $\Gamma$  is a bijection $\varphi\colon \Gamma\to \Gamma$  for which the mapping $x \mapsto x+\varphi(x)$
 is a bijection. In this paper we consider the existence of  a complete mapping $\varphi$ of $\Gamma$ and a partition $S_1,S_2,\ldots S_t$  of elements of $\Gamma$, such that  $\sum_{s\in S_i}s=\sum_{s\in S_i}\varphi(s)=0$ for every $i$, $1 \leq i \leq t$.
 
A $\Gamma$-magic rectangle  set $MRS_{\Gamma}(a, b; c)$  of order $abc$ is a collection of $c$ arrays $(a\times b)$
whose entries are elements of group $\Gamma$ of order $abc$, each appearing once, with all row sums in
every rectangle equal to a constant $\omega\in \Gamma$ and all column sums in every rectangle equal to a
constant $\delta \in \Gamma$.

While a complete characterization of 
MRS$_\Gamma(a,b;c)$ exists for cases where $\{a,b\}\not=\{2k+1,2^{\alpha}\}$, the scenario where $\{a,b\}=\{2k+1,2^{\alpha}\}$ remains unsolved for  $\alpha>1$. Using the partition of $\Gamma$ into zero-sum sets by complete mappings, we give some sufficient conditions that a  $\Gamma$-magic rectangle set MRS$_{\Gamma}(2k+1, 2^{\alpha};c)$ exists.

\noindent\textbf{Keywords:} magic square, magic rectangle  set, complete mapping,  magic constant, Abelian group \\
\noindent\textbf{MSC:} 05B15, 05E99, 05B30
\end{abstract}

 \section{Introduction}

Assume $\Gamma$ is an Abelian group of order $n$ with the operation denoted by $+$.  For convenience
we will write $ka$ to denote $a + a + \ldots + a$ (where the element $a$ appears $k$ times), $-a$ to denote the inverse of $a$ and
we will use $a - b$ instead of $a+(-b)$.  Moreover, the notation $\sum_{a\in S}{a}$ will be used as a short form for $a_1+a_2+a_3+\dots$, where $a_1, a_2, a_3, \dots$ are all elements of the set $S$. The identity element of $\Gamma$ will be denoted by $0$. Recall that an element $\iota\in\Gamma$ of order 2 (i.e., $\iota\neq 0$ and $2\iota=0$) is called an \emph{involution}. 
For convenience, let $\gr$ denote the set consisting of all Abelian groups which are
of odd order or contain more than one involution.
The exponent of a group $\exp(\Gamma)$ is defined as the least common multiple of the orders of all elements of the group. Note that if $\Gamma$ is not cyclic, then  $e(\Gamma)<|\Gamma|$.

A \textit{complete mapping} of a group $\Gamma$ is defined
as $\varphi\in$Bij$(\Gamma)$ (set of all bijections from $\Gamma$ to itself) that the mapping  $\theta\colon g \mapsto  g\varphi(g)$ is
also bijective \cite{Hall}. (Some authors refer to $\theta$, rather than $\varphi$, as the complete
mapping).  
Originally, complete mappings were introduced by Mann in 1942 as a tool for constructing mutually orthogonal Latin squares  \cite{Mann}. For finite Abelian groups, it was proved the following:
\begin{thm}[\cite{Hall,Hall2}]
A finite Abelian $\Gamma$ group has a complete
mapping if and only if $\Gamma\in\gr$.
\end{thm}
Many properties of complete mappings were studied in the literature, see for example \cite{Evans, Tan3}. Let us explain one of these concerning a partition of $\Gamma$ via complete mappings. 
Friedlander,  Gordon Tannenbaum asked whether there exists a complete mapping $\varphi$ of an Abelian group $\Gamma$ which fixes the identity element
and permutes the remaining elements as a product of disjoint $k$-cycles, see~\cite{Tan3}. The problem is still open in general. 

 A subset $S$ of $\Gamma$ is referred to as a \textit{zero-sum subset} if $\sum_{g\in S}g=0$. Note that if the set $\Gamma\setminus\{0\}$ has  a partition $S_1,S_2,\ldots S_t$, such that $|S_i|=3$ and $\sum_{s\in S_i}s=0$ for every  $i=1,2,\ldots,(|\Gamma|-1)3$, then there exists a complete mapping of $\Gamma$ that fixes the identity element, and permutes the remaining elements as products of disjoint cycles of length $3$. By the partition we can write $\Gamma=\{0\}\cup\bigcup_{i=1}^{(|\Gamma|-1)/3}\{x_0^i,x_1^i,x_2^i\},$ where $x_0^i+x_1^i+x_2^i=0$ for any $i=1,2,\ldots,(|\Gamma|-1)/3.$ Set now $\varphi(x_j^i)=x_{j+1}^i$ for $j=0,1,2$, $i=1,2,\ldots,(|\Gamma|-1)/3$, where the subscripts are taken modulo $3$. Note that $\varphi(x_j^i)+x_j^i=x_{j+2}^i+x_{j}^i=-x_{j+3}^i$.
 
This paper considers a partition of $\Gamma$ into zero-sum sets by complete mappings. Note, that $\sum_{g\in\Gamma}g=0$ if and only if $\Gamma\in\gr$. We show that for $\Gamma\in\gr$, there exist $m$ that divides $|\Gamma|$,  a complete mapping $\varphi$ of $\Gamma$ and a partition $S_1,S_2,\ldots S_t$  of elements of $\Gamma$, such that $|S_i|=m$ and $\sum_{s\in S_i}s=\sum_{s\in S_i}\varphi(s)=0$ for every  $1 \leq i \leq t$. We will apply this result for constructing a $\Gamma$-magic rectangle set.

 A \textit{magic square} of order $n$ is an $n\times n$ array with entries $1,2,\ldots,n^2$, each
appearing once, such that the sum of each row, column, and both main diagonals
is equal to $n(n^2+1)/2$.  The earliest known magic square is a $3\times3$ magic square called \textit{Lo Shu magic square} and can be traced in Chinese literature as far back as 2800 B.C. Since then, certainly, many people studied magic squares. For a  survey of magic squares, see Chapter 34 in \cite{ColDin}. Magic rectangles are a natural generalization of magic squares. 
A {\it magic rectangle} MR$(a, b)$ is an $a \times b$ array with entries from the set
$\{1, 2, \ldots , ab\}$, each appearing once, with all its row sums equal to a constant $\delta$ and with
all its column sums equal to a constant $\eta$. The following theorem was proved in \cite{Har1,Har2}:

\begin{thm}[\cite{Har1,Har2}] A magic rectangle MR$(a, b)$ exists if and only if $a, b > 1$, $ab > 4$,
and $a \equiv b \pmod 2$.\label{Har}
\end{thm}

In \cite{Fro}, Froncek introduced the concept of magic rectangle sets. A \textit{magic rectangle set} M= MRS$(a, b; c)$ is a collection of $c$ arrays ($a\times b$)
whose entries are elements of $\{1, 2, \ldots , abc\}$, each appearing once, with all row sums in
every rectangle equal to a constant $\delta=b(abc+1)/2.$ and all column sums in every rectangle equal to a
constant $\eta=a(abc+1)/2.$ It was shown the following.
\begin{thm}[\cite{Fro2}] For $a,b > 1$ and $ab > 4$, a magic rectangle set MRS$(a, b; c)$ exists
if and only if $a, b \equiv 0 \pmod 2$ or $abc \equiv 1 \pmod 2$.\label{thmFro}
\end{thm}

The following generalization magic rectangle set was introduced (\cite{ref_CicIWOCA}).
\begin{dfn}
A $\Gamma$-\textit{magic rectangle  set} MRS$_{\Gamma}(a, b; c)$ on a group $\Gamma$ of order $abc$ is a collection of $c$ arrays $(a\times b)$
whose entries are elements of group $\Gamma$, each appearing once, with all row sums in
every rectangle equal to a constant $\omega\in \Gamma$ and all column sums in every rectangle equal to a
constant $\delta \in \Gamma$. 
\end{dfn}
\begin{thm}[\cite{CicHin}]\label{main} Let   $\{a,b\}\neq\{2^{\alpha},2l+1\}$ for  any natural numbers $\alpha,l>0$. A  $\Gamma$-magic rectangle set MRS$_{\Gamma}(a, b;c)$ exists if and only if $a$ and $b$ are both even or $\Gamma\in\gr$.
\end{thm}

\begin{thm}[\cite{CicHin}]\label{rectangle} Let   $a,b>1$.
A  $\Gamma$-magic rectangle set  MRS$_{\Gamma}(a, b,1)$ exists if and only if $a$ and $b$ are both even or $\Gamma\in\gr$.
\end{thm}

\begin{obs}[\cite{CicHin}]\label{dwa}Let $k,c$ be positive natural numbers and $\Gamma$ be an Abelian group  of order $(4k+2)c$.
There does not exist a $\Gamma$-magic rectangle  set MRS$_{\Gamma}(2k+1, 2;c)$.
\end{obs}
\begin{obs}[\cite{ref_CicIWOCA}]\label{codd}  If $a$ is even, $b$ is odd then for any $c$ and an Abelian group $\Gamma$ having exactly one involution, $|\Gamma|=abc$ there does not exist a $\Gamma$-magic rectangle set MRS$_{\Gamma}(a, b; c)$.
\end{obs}

A complete characterization of MRS$_{\Gamma}(a, b; c)$ is known for $\{a,b\}\neq\{2k+1,2^{\alpha}\}$ where  $k$ and $\alpha$ are the positive integers. However, the case where $\{a,b\}=\{2k+1,2^{\alpha}\}$  remains unsolved for $\alpha>1$ in general. So far only the following is known.

\begin{thm}[\cite{CicHin2}]\label{main2} A necessary and sufficient condition for existence a
MRS$_{\Gamma}(2k + 1, 4; 4l + 2)$ is that the Abelian group $\Gamma$ has more than one involution.
\end{thm}
It was stated:
\begin{conj}[\cite{CicHin}]\label{conjectureST}Let $a,b> 1$.
 A  $\Gamma$-magic rectangle set MRS$_{\Gamma}(a, b;c)$ exists if and only if $a$ and $b$ are both even or $\Gamma\in\gr$ and $\{a,b\}\neq\{2k+1,2\}$.
\end{conj}

One can check, that  in Theorem~\ref{main2}  there has to be $\Gamma\cong A\oplus (\zet_2)^3$ or 
 $\Gamma\cong A\oplus \zet_2\oplus \zet_4$, for $|A|$ odd. Thus $\exp(\Gamma)\equiv 2,4,6\pmod 8$. In this article, we have made progress in proving the above Conjecture~\ref{conjectureST}.
Namely, we generalized Theorem~\ref{main2} for any group $\Gamma$ such that $\exp(\Gamma)\not \equiv 0\pmod 8$. Namely,  we showed that  for any $\Gamma\cong A\oplus (\zet_2)^{\alpha}\oplus (\zet_4)^\beta$ with  $|A|$ odd
a necessary and sufficient condition for the existence a
MRS$_{\Gamma}(2k + 1, 4; (2l + 1)2^{\alpha+\beta-2})$ is that the group $\Gamma$ has more than one involution.

\section{Preliminaries}

 A non-trivial
finite group has elements of order $2$ (involutions) if and only if the order of the group is even. The fundamental theorem of finite Abelian groups states that a finite Abelian
group $\Gamma$ of order $n$ can be expressed as the direct product of cyclic subgroups of prime-power order. This implies that
$$\Gamma\cong\zet_{p_1^{\alpha_1}}\oplus\zet_{p_2^{\alpha_2}}\oplus\ldots\oplus\zet_{p_k^{\alpha_k}}\;\;\; \mathrm{where}\;\;\; n = p_1^{\alpha_1}\cdot p_2^{\alpha_2}\cdot\ldots\cdot p_k^{\alpha_k}$$
and $p_i$ for $i \in \{1, 2,\ldots,k\}$ are not necessarily distinct primes. This product is unique up to the order of the direct product. When $p$ is the number of these cyclic components whose order is a multiple of $2$, then $\Gamma$ has $2^p-1$ involutions. In particular, if $n \equiv 2 \pmod 4$, then $\Gamma\cong \zet_2\oplus \Lambda$ for some
Abelian group $\Lambda$ of odd order $n/2$.  Moreover, every cyclic group of an even order has exactly one involution. \\

If $H$ is a subgroup of $\Gamma$ then we write $H<\Gamma$. By $\langle g \rangle$ we denote the subgroup generated by $g$ in the group $\Gamma$.


Since the properties and results in this paper are invariant under an isomorphism between groups, we only need to consider one group in each isomorphism class.

\section{Complete mappings}
It was proved the following.
\begin{thm}[\cite{CicZ}]\label{mZSP} Let $\Gamma$  be an Abelian group of order $n$ such that  $\Gamma\in\gr$. Let $m>1$ divide $n$.
There exists a partition $S_1,S_2,\ldots S_t$  of elements of $\Gamma$, such that $|S_i|=m$ and $\sum_{s\in S_i}s=0$ for every  $i=1,2,\ldots,t$.
\end{thm}
We will start with the groups such that $|\Gamma|$ odd.
\begin{lem}\label{partitionn}
Let $\Gamma$  be an Abelian group such that   $|\Gamma|=n$ is odd. Let $m>1$ divide $n$.
There exist a complete mapping $\varphi$ of $\Gamma$ and a partition $S_1,S_2,\ldots S_t$  of elements of $\Gamma$, such that $|S_i|=m$ and $\sum_{s\in S_i}s=\sum_{s\in S_i}\varphi(s)=0$ for every  $1 \leq i \leq t$.
\end{lem}
\textit{Proof.} By Theorem~\ref{mZSP} there exists a partition $S_1,S_2,\ldots S_t$  of elements of $\Gamma$, such that $|S_i|=m$ and $\sum_{s\in S_i}s=0$ for every  $i=1,2,\ldots,t$. Let $\varphi(g)=g$ for any $g\in \Gamma$, then $\varphi+id\mapsto 2g$, which is an automorphism since $|\Gamma|$ is odd.~\qed

For $\Gamma\in \gr$ such that $|\Gamma|=2^n$  using the idea from~\cite{CicZ} and \cite{Zeng} we will show the following:
\begin{thm}\label{partition}
Let $\Gamma$ of order $2^n$, $n>1$ be an Abelian group such that  $\Gamma\in\gr$. 
Let $$m=\begin{cases}2\exp(\Gamma)& \text{if}\;\;\exp(\Gamma)=\frac{|\Gamma|}{2},\\
\max\{4,\exp(\Gamma)\}& \text{if}\;\; \exp(\Gamma)<\frac{|\Gamma|}{2}.\\
\end{cases}$$
There exist a complete mapping $\varphi$ of $\Gamma$ and a partition $S_1,S_2,\ldots S_t$  of elements of $\Gamma$, such that $|S_i|=m$ and $\sum_{s\in S_i}s=\sum_{s\in S_i}\varphi(s)=0$ for every  $1 \leq i \leq t$.
\end{thm}
\textit{Proof.} For $\Gamma$ such that $exp(\Gamma)=\frac{|\Gamma|}{2}$, there is $\Gamma\cong \zet_{2^\beta}\oplus\zet_{2}$ and the conclusion is obvious, since $\sum_{g\in \Gamma}g=0$ and $\exp(\Gamma)=2^\beta$.

From now we assume that $\Gamma\not\cong \zet_{2^\beta}\oplus\zet_{2}$ for $\beta\geq 1$.

The proof is by induction on $|\Gamma|$. We deal with some base cases. Suppose first that $\Gamma\cong(\zet_2)^m$ for some $m>2$. Take a subgroup $\Gamma_0\cong\zet_2\oplus\zet_2$ of $\Gamma$, then  there exists a complete mapping $\varphi_0$ of $\Gamma_0$ and the partition consists of only one set $S=\Gamma_0$. Let $R=\{\iota_1,\iota_2,\ldots,\iota_{2^{m-2}}\}$, be the coset representative of $\Gamma/\Gamma_0$. Note that for any of that representative $ x\in R$ there is $2x=0\in\Gamma$. Let $\Gamma\ni a=a_0+x$ for $x\in R$ and $s_0\in\Gamma_0$, set $\varphi(a)=\varphi_0(a_0)+x$. The partition of $S_1,S_2,\ldots,S_{2^{m-2}}$ defined as $S_j=\{a_0+\iota_j,\colon a_0\in\Gamma_0\}$ is the desired partition of $\Gamma$.

Let now $$\Gamma \cong \zet_{2^{\alpha_1}}\oplus\zet_{2^{\alpha_2}}\oplus\ldots\oplus\zet_{2^{\alpha_k}}$$
  with $\alpha_1\geq2$ and $\alpha_2+\alpha_3\geq2$, $\alpha_1\geq\alpha_2\ldots\geq\alpha_k$. Then there exists a subgroup $$\Gamma_0\cong\langle (2^{\alpha_1-1},2^{\alpha_2-1},1,1,\ldots,1)\rangle\cong \zet_{2^{\alpha_1-1}}\oplus\zet_{2^{\alpha_2-1}}\oplus\zet_{2^{\alpha_3}}\ldots\oplus\zet_{2^{\alpha_k}} \in\gr$$ of $\Gamma$   such $\Gamma/\Gamma_0\cong \zet_2\oplus\zet_2$.  
  Observe that $\exp(\Gamma)=2\exp(\Gamma_0)$.   Let
  $$m_0=\begin{cases}2\exp(\Gamma_0)& \text{if}\;\;\exp(\Gamma_0)=\frac{|\Gamma_0|}{2},\\
\max\{4,\exp(\Gamma_0)\}& \text{if}\;\; \exp(\Gamma_0)<\frac{|\Gamma_0|}{2}.\\
\end{cases}$$
  
   By the induction hypothesis, there is a complete mapping  $\varphi_0\in$Bij$(\Gamma_0)$ and a partition $S_1',S_2',\ldots S_t'$  of elements of $\Gamma_0$, such that $|S_i'|=m_0$ and $\sum_{s\in S_i'}s=\sum_{s\in S_i}\varphi_0(s)=0\in\Gamma$ for every $i=1,2,\ldots,t$. 

 Choose a set of coset representatives for $\Gamma/\Gamma_0$ to be $\{0,c,d,-c-d\}$. Note that
$$(c+\Gamma_0)\bigcup(d+\Gamma_0)\bigcup(-c-d+\Gamma_0)=\bigcup_{b\in\Gamma_0}\{c+b,d+\phi_0(b),-c-d+\varphi_0(b)\},$$
and every subset $\{c+b,d+\phi_0(b),-c-d+\varphi_0(b)\}$ is zero-sum. 

Let us define the partition as follows:
$S_i^0=\{g\colon g\in S_i'\}$, $S_i^1=\{g+c\colon g\in S_i'\}$, $S_i^2=\{g+d\colon g\in S_i\}$ and $S_i^3=\{g-c-d\colon g\in S_i\}$ for $i=1,2,\ldots,t'$.
Note that  $\sum_{a\in S_i^0}a=\sum_{g\in S_i}\varphi_0(g)=0\in\Gamma$, whereas for $j\neq 0$ we have $\sum_{a\in S_i^j}a=\sum_{g\in S_i}\varphi_0(g)+m_0e=0+m_0e\in\Gamma$ for some  $e\in\{c,d,-c-d\}$. Moreover, if $e\in S_i^1$, then $e-c+d\in S_i^2$ and $e-2c-d\in S_i^3$. 
Define now $\varphi$ as follows:
$$\varphi(a)=\left\{\begin{array}{ccl}
\varphi_0(a) & \mathrm{if} & a\in S_0, \\
\varphi_0(a)-c+d & \mathrm{if} & a\in S_i^1, \\
\varphi_0(a)-2d-c & \mathrm{if} & a\in S_i^2, \\
\varphi_0(a)+d+2c & \mathrm{if} & a\in S_i^3 \end{array}\right.$$

Note that $$\varphi(a)+a=\left\{\begin{array}{ccl}
\varphi_0(a) +a & \mathrm{if} & a\in S_0, \\
\varphi_0(a)+a+d & \mathrm{if} & a\in S_i^1, \\
\varphi_0(a)+a-d-c & \mathrm{if} & a\in S_i^2, \\
\varphi_0(a)+a+c & \mathrm{if} & a\in S_i^3 \end{array}\right.$$
Thus $\varphi+id\in$Bij$(\Gamma)$. Moreover, observe that for $j=1,2$ there is $\sum_{a\in S_i^j}\varphi(a)=\sum_{a\in S_i^{j+1}}a$ and  $\sum_{a\in S_i^3}\varphi(a)=\sum_{a\in S_i^{1}}a$.
Thus for $j\neq 0$, $i=1,2,\ldots,t'$  $\sum_{a\in S_i^j}\varphi(a)=0+m_0e\in\Gamma$ for some  $e\in\{c,d,-c-d\}$.

Observe that since $t'= |\Gamma|/16$, $t'$ is odd if and only if $\Gamma\cong \zet_4\oplus\zet_2\oplus\zet_2$, so $m=m_0$. Therefore for  $m_0=m/2$ we can assume that $t'$ is even. Moreover  $t=t'$ if $m_0=m$ and $t=t'/2$ if  $m_0=m/2$. Set a partition
$$S_i=\begin{cases}S_{\lfloor{i/t'}\rfloor}^{i\pmod 4},&\;\text{if}\;\; m_0=m,\\
S_{\lfloor{2i/t'}\rfloor}^{i\pmod 4}\cup S_{\lfloor{(2i+1)/t'}\rfloor}^{i\pmod 4}&\;\text{if}\;\; m_0=m/2.
\end{cases}$$
for  $i=1,2,\ldots,t$.
One can easily see that  $\sum_{a\in S_i}\varphi(a)=0+me=0\in\Gamma$ for $e\in\{0,c,d,-c-d\}$ and $i=1,2,\ldots,t$. This finishes the proof. ~\qed

\begin{thm}\label{partitiong}
Let $\Gamma \cong L\oplus H\in\gr$ be such that $|H|$ is odd, $|L|=2^{\eta}$ for some natural number $\eta$. Let $l=2\exp(L)$ for $\exp(\Gamma)=|\Gamma|/2$ and $l=\exp(L)$ otherwise. Let $k>1$, $k||H|$ and $m=kl$. There exist a complete mapping $\varphi$ of $\Gamma$ and a partition $S_1,S_2,\ldots S_t$  of elements of $\Gamma$, such that $|S_i|=m$ and $\sum_{s\in S_i}s=\sum_{s\in S_i}\varphi(s)=0$ for every  $1 \leq i \leq t$.
\end{thm}
\textit{Proof.} Let $S_1^1,S_2^1,\ldots S_{t_1}^1$ be partition of $L$, such that $|S_i^1|=l$ and $\sum_{s\in S_i^1}s=\sum_{s\in S_i^1}\varphi_1(s)=0\in L$ for every $i$, $1 \leq i \leq t_1$, which exists by Lemma~\ref{partition}. Let $S_1^2,S_2^2,\ldots S_{t_2}^2$ be a partition of $H$, such that $|S_i^2|=k$ and $\sum_{s\in S_i^2}s=\sum_{s\in S_i^2}\varphi_2(s)=0\in L$ for every $i$, $1 \leq i \leq t_2$, which exists by Lemma~\ref{partitionn}.

Let the complete mapping for $\Gamma_1$ be $\varphi_1$ and, whereas  for $\Gamma_2$ be $\varphi_2$. Then set

$$\varphi=(\varphi_1,\varphi_2): \Gamma_1\oplus \Gamma_2\rightarrow\Gamma_1\oplus \Gamma_2,\;\; (a_1,a_2)\mapsto (\varphi_1(a_1),\varphi_2(a_2)).$$

Then the desired partition of $\Gamma$ contains the Cartesian products of $S_i^1\times S_j^2$, $i=1,2,\ldots,t_1$, $i=j,2,\ldots,t_2$.~\qed

As a corollary, we obtain the following immediately:
\begin{cor}
Let $\Gamma\in \gr$. If there exist  prime numbers $p$ and $q$,  such that $pq$ divides $|H|$ or $\exp(\Gamma)<|\Gamma|/2$. Then there exist a positive integer $m<|\Gamma|$, a complete mapping $\varphi$ of $\Gamma$ and a partition $S_1,S_2,\ldots S_t$  of elements of $\Gamma$, such that $|S_i|=m$ and $\sum_{s\in S_i}s=\sum_{s\in S_i}\varphi(s)=0$ for every  $1 \leq i \leq t$.
\end{cor}

In \cite{Wal}, Marr and Wallis defined a Kotzig array as a $j\times k$ grid where each row is a permutation of $\{1,2,\ldots,k\}$ and each column has the same sum $j(k+1)/2$. Notably, a Latin square is a specific type of Kotzig array.
\begin{thm}[\cite{Wal}]\label{Wal}
A Kotzig array of size $j \times k$ exists if and only if
$j>1$ and $j(k - 1)$ is even. 
\end{thm}

In \cite{CicZ}, the author introduced a generalization of Kotzig arrays and established necessary and sufficient conditions for their existence. Specifically, for an Abelian group $\Gamma$ of order $k$, a $\Gamma$-Kotzig array of size $j\times k$ is defined as a $j\times k$ grid where each row is a permutation of elements of $\Gamma$ and each column has the same sum. Moreover, without loss of generality, it can be assumed that all elements in the first column are 
$0$ (thus, all column sums are $0$).
\begin{thm}[\cite{CicZ}]\label{Kotzig}
A $\Gamma$-Kotzig array of size $j \times k$ exists if and only if
$j>1$ and $j$ is even or $\Gamma\in \gr$. 
\end{thm}
In this paper we introduce a $\Gamma$-Kotzig arrays  set KAS$_{\Gamma}(j, m; k/m)$ a collection of $k/m$ arrays $j\times m$ that is made from a partition of a $\Gamma$-Kotzig array of size $j\times k$ into $k/m$ arrays of $j\times m$ such that the sum of all elements in rows and columns is $0$.

We will start with some obvious observations.
\begin{obs}\label{glue}
Let $\Gamma\in\gr$ be of order $n=abc$ and $a>1$.
If there exists a $\Gamma$-Kotzig arrays  set KAS$_{\Gamma}(j, a; bc)$, then there exists a $\Gamma$-Kotzig arrays  set KAS$_{\Gamma}(j, ab; c)$.
\end{obs}
\textit{Proof.} There exists a $\Gamma$-Kotzig arrays  set KAS$_{\Gamma}(j, a; bc)$. To construct
one of KAS$_{\Gamma}(j, ab; c)$, we simply take $b$ of  KAS$_{\Gamma}(j, a; bc)$ arrays and ”glue” them into an array  $j\times ab$.~\qed

\begin{lem}
Let $\Gamma \in\gr$ be of order $n$ and $m>1$ divide $n$.
There exists a $\Gamma$-Kotzig arrays  set KAS$_{\Gamma}(2, m; n/m)$.
\end{lem}
\textit{Proof.} There exists a partition $S_1,S_2,\ldots S_{n/m}$  of elements of $\Gamma$, such that $|S_i|=m$ and $\sum_{s\in S_i}s=0$ for every  $1 \leq i \leq n/m$
by Theorem~\ref{mZSP}.

Let $S_i=\{g_1^i,g_2^i,\ldots,g_m^i\}$ for $i=1,2,\ldots,n/m$.
Denote by $k^s_{i,j}$ the entry in the $i$-th row and $j$-th column of the $s$-th array in the  $\Gamma$-Kotzig arrays  set KAS$_{\Gamma}(2, m; n/m)$. 
Let $x_{1,i}^s=g_i^s$, $x_{2,i}^s=-g_i^s$ for $i=1,2,\ldots,m$, $s=1,2,\ldots,n/m$.~\qed

The below lemmas will be very useful for the result on magic rectangle sets:

\begin{lem}
Let $\Gamma \in\gr$ be of order $n=2^{\alpha}$ and $$m=\begin{cases}2\exp(\Gamma)& \text{if}\;\;\exp(\Gamma)=\frac{|\Gamma|}{2},\\
\max\{4,\exp(\Gamma)\}& \text{if}\;\; \exp(\Gamma)<\frac{|\Gamma|}{2}.\\
\end{cases}$$

There exists a  a $\Gamma$-Kotzig arrays  set KAS$_{\Gamma}(3, m;n/m)$. 
\end{lem}
\textit{Proof.} 
 Let $t=n/m$. By Lemma~\ref{partition} there exist a complete mapping $\varphi$ of $\Gamma$ and a partition $S_1,S_2,\ldots S_t$  of elements of $\Gamma$,  such that $|S_i|=m$ and $\sum_{s\in S_i}s=\sum_{s\in S_i}\varphi(s)=0$ for every  $1 \leq i \leq t$.

Let $S_i=\{g_1^i,g_2^i,\ldots,g_m^i\}$ for $i=1,2,\ldots,t$.
Denote by $k^s_{i,j}$ the entry in the $i$-th row and $j$-th column of the $s$-th array in the  $\Gamma$-Kotzig arrays  set KAS$_{\Gamma}(2, m; t)$. 
Let $x_{1,i}^s=g_i^s$, $x_{2,i}^s=\varphi(g_i^s)$ and $x_{3,i}^s=-g_i^s-\varphi(g_i^s)$ for $i=1,2,\ldots,m$, $s=1,2
,\ldots,t$.~\qed

\begin{thm}\label{KAS}
Let $j>1$, $\Gamma \in\gr$ be of order $n=2^{\alpha}$ and $$m=\begin{cases}2\exp(\Gamma)& \text{if}\;\;\exp(\Gamma)=\frac{|\Gamma|}{2},\\
\max\{4,\exp(\Gamma)\}& \text{if}\;\; \exp(\Gamma)<\frac{|\Gamma|}{2}.\\
\end{cases}$$
There exists a $\Gamma$-Kotzig arrays  set KAS$_{\Gamma}(j, m; n/m)$.  
\end{thm}
\textit{Proof}.  Assume first $j$ is even. To construct
a $\Gamma$-Kotzig array of size $j\times k$, we simply take $j/2$ set the KAS$_{\Gamma}(2, m; n/m)$ of size $2\times k$ and "glue" them into the KAS$_{\Gamma}(j, m; n/m)$ of size.

To construct a  KAS$_{\Gamma}(j, m; n/m)$ for $j$ odd, we simply take $(j-3)/2$ of KAS$_{\Gamma}(2, m; n/m)$, one KAS$_{\Gamma}(3, m; n/m)$ and "glue" them into a KAS$_{\Gamma}(j, m; n/m)$.~\qed


\section{Main results}

The following two lemmas are also useful tools in the proof of the main result.

\begin{lem}[\cite{CicHin}]\label{lemgl}
Let $\Gamma$ be a group of order $abc_1c_2$. Let  $\Gamma\cong \Gamma_0\times H$ for some group $H\in \gr$ of order $c_2$. If there exists a $\Gamma_0$-magic rectangle set  MRS$_{\Gamma_0}(a, b;c_1)$   with the  column sum $\delta_0$ and the row sum $\omega_0$, then there exists a $\Gamma$-magic rectangle set MRS$_{\Gamma}(a, b; c_1c_2)$ with the  column sum $(\delta_0,0)$ and the row sum $(\omega_0,0)$.
\end{lem}

\begin{lem}[\cite{CicHin}]\label{lemgl2}
Let $\Gamma\cong \gA\oplus \zet_{2k+1}$ for some Abelian group $\gA$. Let $h$ be a natural number that divides $2k+1$ and $\Gamma_0\cong\gA\times\left\langle  h\right\rangle$. If there exists a $\Gamma_0$-magic rectangle set  MRS$_{\Gamma_0}(a, b;c_1)$, then there exists a $\Gamma$-magic rectangle set MRS$_{\Gamma}(a, b; c_1h)$.
\end{lem}

We show the following:

\begin{lem}\label{p3} 
Let $\Gamma\cong \zet_3\oplus \Delta$  with  $|\Delta|=2^{\alpha}$ and $exp(\Gamma)\leq 4$, $\Delta\in \gr$. Then there exists a MRS$_{\Gamma}(3, 4, 2^{\alpha-3})$.
\end{lem}
\textit{Proof.} The proof is by induction on $|\Delta|$.
Let $\Gamma_1\cong \zet_3\oplus\zet_2\oplus\zet_2$, $\Gamma_2\cong\zet_3\oplus\zet_4\oplus\zet_2$ and $\Gamma_3\cong \zet_3\oplus\zet_2\oplus\zet_2\oplus\zet_2$.
We deal with some base cases. If $|\Delta|\in\{4,8\}$, then by Theorems~\ref{rectangle} and \ref{main2}  there exist MRS$_{\Gamma_1}(3, 4, 1)$,  MRS$_{\Gamma_2}(3, 4, 2)$ and 
MRS$_{\Gamma_3}(3, 4, 2)$. Assume now that $|\Delta|\geq 16$. Then there exists  a subgroup $H$ of $\Delta$ such that  $\Gamma/H\cong \zet_2\oplus\zet_2$. By induction there exists a MRS$_{H}(3, 4, 2^{n-4})$ with  row sum  $\omega\in \Gamma$ and column sum $\delta\in\Gamma$.  Choose a set of coset representatives for $\Gamma/H$ to be $\{2(c+d),c,d,c+d\}$. Let $c_2=c$, $c_3=d$ and $c_4=c+d$.
 Denote by $y_{i,j}^s$ the entry in the $i$-th row and $j$-th column of the $s$-th rectangle in the $H$-magic rectangle set MRS$_{\Gamma}(3, 4; 2^{n-4})$. and by $x^s_{i,j}$ the entry in the $i$-th row and $j$-th column of the $s$-th rectangle in the $\Gamma$-magic rectangle set MRS$_{\Gamma}(3, 4; 2^{n-2})$.   Let $x_{i,j}^s= y_{i,j}^{\left\lfloor s/4\right\rfloor}+2(c+d)$ for $s\equiv 1\pmod 4$, $i=1,2,3$, $j=1,2,3,4$. Note that for that kind of rectangles the row sum is $\omega+4\cdot2(c+d)=\omega$ and column sum $\delta+6\cdot2(c+d)=\delta+2(c+d)$.
 
 Let
  $$x_{i,j}^{s}=\begin{cases}(y_{i,j}^{\left\lfloor s/4\right\rfloor},c)&i=1,j=1,2,3,4,\\
 (y_{i,j}^{\left\lfloor s/4\right\rfloor},d)&i=2,j=1,2,3,4,\\
(y_{i,j}^{\left\lfloor s/4\right\rfloor},c+d)&i=3,j=1,2,3,4,\\\end{cases}$$
for  $s\equiv 2\pmod 4$,

 $$x_{i,j}^{s}=\begin{cases}(y_{i,j}^{\left\lfloor s/4\right\rfloor},d)&i=1,j=1,2,3,4,\\
 (y_{i,j}^{\left\lfloor s/4\right\rfloor},d+c)&i=2,j=1,2,3,4,\\
(y_{i,j}^{\left\lfloor s/4\right\rfloor},c)&i=3,j=1,2,3,4,\\\end{cases}$$
 for $s\equiv 3\pmod 4$ and 
 $$x_{i,j}^{s}=\begin{cases}(y_{i,j}^{\left\lfloor s/4\right\rfloor},c+d)&i=1,j=1,2,3,4,\\
 (y_{i,j}^{\left\lfloor s/4\right\rfloor},c)&i=2,j=1,2,3,4,\\
(y_{i,j}^{\left\lfloor s/4\right\rfloor},d)&i=3,j=1,2,3,4,\\\end{cases}$$
for $s\equiv 0\pmod 4$.
 Note that for that kind of rectangles the row sum is $\omega$  and column sum $\delta+2(c+d)$.~\qed

\begin{lem}\label{pnotexp} 
Let $p>3$ be a prime and $\Gamma\cong \zet_p\oplus \Delta$  with  $|\Delta|=2^{\alpha}=n$ and $\exp(\Delta)<|\Delta|/2$. There exists a MRS$_{\Gamma}(p, \exp(\Delta);n/\exp(\Delta))$, if $\gcd(p-1,\exp(\Delta))=1$ and a MRS$_{\Gamma}(p, 2\exp(\Delta);n/(2\exp(\Delta)))$ otherwise.
\end{lem}
\textit{Proof.}
Note that if $\gcd(p-1,\exp(\Delta))=p$, then  $\gcd(p-1,2\exp(\Delta))=1$. Let $m=\exp(\Delta)$ if $\gcd(p-1,\exp(\Delta))=1$ and $m=2\exp(\Delta)$ otherwise. Let $f(x)=-(m-1)x$ for any $x\in \zet_p$. Note that since $\gcd(m-1,p)=1$ the mapping $f$ is an automorphism.
Let $f(x)=-(m-1)x$ for any $x\in \zet_p$.  By Theorem~\ref{KAS} and Observation~\ref{glue} there 
exists a   $\Delta$-Kotzig arrays  set KAS$_{\Delta}(p, m; n/m)$. Denote by $k^s_{i,j}$ the entry in the $i$-th row and $j$-th column of the $s$-th array in the  $\Gamma$-Kotzig arrays  set KAS$_{\Gamma}(p, m; n/m)$.

Denote  by  $x^s_{i,j}$ the entry in the $i$-th row and $j$-th column of the $s$-th rectangle in the $\Gamma_0$-magic rectangle set MRS$_{\Gamma}(p, m; n/m)$.

 Set
$$x_{i,j}^s=\left\{
\begin{array}{lcl}
(f(i-1),k_{i,j}^{\lfloor\frac{ms}{n}\rfloor}), & \text{if} & j=1, \\
(i-1,k_{i,j}^{\lfloor\frac{ms}{n}\rfloor}), & \text{if} & j\in\{2,3,\ldots,m\}, \\
\end{array}%
\right.
$$
for $i=1,2,\ldots,p$.

In this  case  $x_{i,j}^s\neq x_{i',j'}^{s'}$ for $(i,j,s)\neq (i',j',s')$ and every column sum is $(\sum_{x\in \zet_p}x,\sum_{i=1}^pk_{i,j})=(0,0)$ and every row sum  is $(0,0)$.~\qed\\

We immediately obtain the following observations:
\begin{obs}\label{pmrs} 
Let $p>2$ be a prime and $\Gamma\cong \zet_p\oplus \Delta$  with  $|\Delta|=2^{\alpha}$ and  $\Delta\in \gr$. Then there exists a MRS$_{\Gamma}(p, 4, 2^{\alpha-2})$.
\end{obs}
\textit{Proof.} If $p=3$ then we apply Lemma~\ref{p3}. Assume that $p>3$. If now  $\exp(\Delta)=|\Delta|/2$, then $\Delta\cong\zet_4\oplus\zet_2$ and we are done by Theorem~\ref{main2}.  For other cases we use Lemma~\ref{pnotexp}, since $\gcd(p,3)=1$.~\qed

\begin{obs}\label{exp} 
Let $p>3$ be a prime and $\Gamma\cong \zet_p\oplus \Delta$  with  $|\Delta|=2^{\alpha}$ and  $\Delta\in \gr$. There exists a MRS$_{\Gamma}(p, 2\exp(\Delta), |\Delta|/(2\exp(\Delta))$.
\end{obs}
\textit{Proof.} If $\exp(\Gamma)=|\Gamma|/2$, then there exists a a MRS$_{\Gamma}(p, 2\exp(\Gamma), 1)$ by Theorem~\ref{rectangle}. Assume that $\exp(\Gamma)<|\Gamma|/2$. We are done by Lemma~\ref{pnotexp} now.~\qed\\

We will finish that section with our main result:
\begin{thm}\label{lem2p} Let $k>1$ and $c$ be odd integers and $\Gamma$ with  $\exp(\Gamma)\not \equiv 0\pmod 8$.  A  $\Gamma$-magic rectangle set MRS$_{\Gamma}(a, b;c)$ exists if and only if $a$ and $b$ are both even or $\Gamma\in\gr$ and $\{a,b\}\neq\{2k+1,2\}$.
\end{thm}
\textit{Proof.}  For $a\equiv b \pmod 2$ or $a$ odd and $b\neq 2^{\alpha}$ we are done by Theorem~\ref{main}. Thus we can assume that $a>1$ is odd and $b=2^\alpha$ for some positive integer $\alpha$. If now $\alpha=1$ or $\Gamma\not\in \gr$ then we are done by Observation~\ref{dwa} or \ref{codd}, respectively. Thus we can assume that $a>1$ and $b=2^\alpha$, $\alpha>2$ and $\Gamma\cong  A \oplus \Delta$ with  $|A|=a$ odd,  $|\Delta|=2^{\alpha}$ and $\Delta\in\gr$, $\exp(\Delta)\leq 4$. Let $a=kc$ for $k>1$, $c$ odd. Note that it is enough to show that 
there exists a MRS$_{\Gamma}(k, 4;2^{\alpha-2}c )$.
By the Fundamental Theorem of Finite Abelian Groups the group $A$  can be written as the direct product: %
$$A\cong\zet_{q_1^{\alpha_1}}\times\ldots\times\zet_{q_w^{\alpha_w}}$$
where $q_1,\ldots,q_w$ are odd (not necessary distinct) primes  and $ck =  q_1^{\alpha_1}\cdot\ldots\cdot q_w^{\alpha_w}$. Without loss of generality, we can assume that  $(2k+1)=q_1\cdot q$ for some odd $q\geq1$.

The strategy for the proof involves using a MRS$_{\zet_{q_1}\oplus\Delta}(q_1, 4;2^{\alpha-2})$ as a "starter" and then "blowing it up" into a MRS$_{\zet_{q_1^{\alpha_1}}\oplus\Delta}(k,4;2^{\alpha-2}c)$ by applying Lemmas~\ref{lemgl} and~\ref{lemgl2}.

Specifically, by Observation~\ref{pmrs} there exists a MRS$_{\zet_{q_1}\oplus\Delta}(q_1, 4;2^{\alpha-2})$. By applying now Lemma~\ref{lemgl2}  we can establish the existence of a MRS$_{\zet_{q_1^{\alpha_1}}\oplus\Delta}(q_1, 4;2^{\alpha-2}q_1^{\alpha_1-1})$ (it can happen that $\alpha_1=1$).  Lemma~\ref{lemgl} further implies the existence of a MRS$_{A\oplus\Delta}(q_1, 4;2^{\alpha-2}kc/q_1)$. To construct  one of  the rectangles from MRS$_{A\oplus\Delta}(k, 4;2^{\beta-2}c)$ we simply take  $k/q_1$ of  MRS$_{A\oplus\Delta^{\beta}}(q_1, 4;2^{\beta-2}kc/q_1)$ rectangles and ''glue'' them together to form a rectangle of size $k \times 4$.~\qed

Note that using a similar method as above and Observation~\ref{exp} we obtain the following:

\begin{obs}\label{exp} 
Let  $\Gamma\cong A\oplus \Delta$  with $|A|$ odd,  $|\Delta|=2^{\alpha}$ and  $\Delta\in \gr$.  If $k>1$ divides, then there exists a MRS$_{\Gamma}(k, 2\exp(\Gamma), |\Gamma|/(2k\exp(\Gamma))$.
\end{obs}

\section{Final Remarks}\label{sec:final}
We have achieved important progress towards proving Conjecture~\ref{conjectureST} by showing that it is true for any $\Gamma\in\gr$ such that $\exp(\Gamma)\not\equiv0\pmod 8$. 
 This proof's crucial idea was a zero-sum group partition via complete mapping. Note that, if for an Abelian group $\Gamma$ of odd order, there is a zero-sum partition into test $S_1,S_2,\ldots,S_t$, then there exists a perfect mapping $\varphi$ of $\Gamma$ such that $\sum_{a\in S_i}\varphi(a)=0\in\Gamma$, namely $\varphi=id_{\Gamma}$. Therefore for an Abelian group $\Gamma$ of an even order we state the open problem:\\
 
\textbf{Open problem:} Let $\Gamma$  be an Abelian group such that   $\Gamma\in \gr$. Determine the smallest possible $m$ such that there exist a complete mapping $\varphi$ of $\Gamma$ and a partition $S_1,S_2,\ldots S_t$  of elements of $\Gamma$, that $|S_i|=m$ and $\sum_{s\in S_i}s=\sum_{s\in S_i}\varphi(s)=0$ for every $i$, $1 \leq i \leq t$.

\section{Statements and Declarations}
This work was  supported by program ''Excellence initiative – research university'' for the AGH University.

\bibliographystyle{plain}

\end{document}